\documentclass{amsart}

\newtheorem{theorem}{Theorem}[section]
\newtheorem{lemma}[theorem]{Lemma}
\newtheorem{proposition}[theorem]{Proposition}
\newtheorem{problem}[theorem]{Problem}

\newtheorem{fact}[theorem]{Fact}

\theoremstyle{definition}

\newtheorem{example}[theorem]{Example}
\usepackage{amsmath,amssymb,amsfonts}

\theoremstyle{remark}
\newtheorem{remark}[theorem]{Remark}

\numberwithin{equation}{section}

\renewcommand{\span}{\mathrm{span}}
\newcommand{\diam}{\mathrm{diam}} 
\newcommand{\dist}{\mathrm{dist}}
\newcommand{\dens}{\mathrm{dens}}
\newcommand{\conv}{\mathrm{conv}}

\newcommand{\R}{\mathbb{R}}

\newcommand{\X}{\mathrm{X}}

\newcommand{\Y}{\mathrm{Y}}
\newcommand{\Z}{\mathrm{Z}}
\newcommand{\B}{\mathbf{B}}
\newcommand{\I}{\mathbf{I}}

\newcommand{\cf}{\mathrm{cf}}

\renewcommand{\S}{\mathbf{S}}

\renewcommand{\ker}{\mathrm{Ker}}
\newcommand{\quotient}{/}
\newcommand{\1}{\boldsymbol{1}}

\begin{document}

\title{Lindel\"{o}f type of generalization of separability in Banach spaces}

\author{Jarno Talponen}
\address{University of Helsinki, Department of Mathematics and Statistics, Box 68, (Gustaf H\"{a}llstr\"{o}minkatu 2b) FI-00014 University 
of Helsinki, Finland}
\email{talponen@cc.helsinki.fi}

\date{\today}

\begin{abstract}
We will introduce the countable separation property (CSP) of Banach spaces $\X$, which is defined as follows: 
For each set $\mathcal{F}\subset \X^{\ast}$ such that $\overline{\span}^{\omega^{\ast}}(\mathcal{F})=\X^{\ast}$ 
there is a countable subset $\mathcal{F}_{0}\subset\mathcal{F}$ with $\overline{\span}^{\omega^{\ast}}(\mathcal{F}_{0})=\X^{\ast}$. 
All separable Banach spaces have CSP and plenty of examples of non-separable CSP spaces are provided. Connections of CSP
with Marku\v{c}evi\v{c}-bases, Corson property and related geometric issues are discussed.
\end{abstract}
\maketitle

\section{Introduction}
In this article we study a class of relatively small non-separable Banach spaces.
If $\X$ is a Banach space and $\X^{\ast}$ is its dual, then a subset $\mathcal{F}\subset\X^{\ast}$ is said to \emph{separate}
$\X$ if for all $x\in\X\setminus\{0\}$ there is $f\in\mathcal{F}$ such that $f(x)\neq 0$. We investigate 
the Banach spaces $\X$, which have the following \emph{countable separation property} (CSP):
\[\mathrm{Whenever}\ \mathcal{F}\subset\X^{\ast}\ \mathrm{separates}\ \X\mathrm{,\ there\ exists\ a\ countable\ separating\ subset}\ 
\mathcal{F}_{0}\subset\mathcal{F}.\]

It turns out that all separable Banach space have CSP but there exist also numerous examples
of non-separable CSP spaces (see the last section for summary). The definition of the countable separation property resembles 
that of the Lindel\"{o}f property. In fact it turns out that if the unit sphere $\S_{\X}$ of a Banach space $\X$ is weakly Lindel\"{o}f, 
then $\X$ has CSP. We will point out various connections between CSP and the established theory of non-separable spaces.
For example, if $\X$ is weakly compactly generated and has CSP, then it follows that $\X$ is already separable.

This study is also closely related to \cite{Gr}, in which certain generalizations of
separability, the so called \emph{Kunen-Shelah properties}, are summarized and discussed. The chain of implications 
\[\mathrm{separability}\implies\mathrm{KS}_{7}\implies\mathrm{KS}_{6}\implies\dots\implies\mathrm{KS}_{0} \] 
was proved there. The only known examples of non-separable spaces having any of the properties $\mathrm{KS}_{7}-\mathrm{KS}_{2}$ 
are constructed under set theoretic axioms extraneous to ZFC. To mention in this context the most important properties from the list, 
$\mathrm{KS}_{4}$, $\mathrm{KS}_{2}$ or $\mathrm{KS}_{1}$ of $\X$ means that there is \emph{no} closed subspace $\Y\subset\X$ admitting an 
$\omega_{1}$-polyhedron, an uncountable biorthogonal system or an uncountable M-basis, respectively, in $\Y$. It turns out that 
$\mathrm{KS}_{4}\implies$CSP$\implies\mathrm{KS}_{1}$. 

In order to motivate the introduction of CSP concept, let us discuss the Kunen-Shelah properties a little further.  
One could ask whether some of the Kunen-Shelah properties of Banach spaces are in fact equivalent to separability. 
Recently Todor\v{c}evi\'{c} \cite{To} answered some old questions about the existence of bases in non-separable Banach spaces. 
In particular, he proved that it is consistent that each non-separable Banach space admits an uncountable biorthogonal
system and hence it is consistent that $\mathrm{KS}_{2}\implies$CSP. 
Actually, by combining the results in \cite[p.1086-1099]{Ne} and \cite{To} one can see that the question about the equivalence 
of separability with $\mathrm{KS}_{2}$ is independent of ZFC.
This positions the class of CSP spaces interestingly with respect to the Kunen-Shelah properties. 
On one hand, CSP is a class of Banach spaces close to $\mathrm{KS}_{2}$ and still containing absolute non-separable representatives. 
On the other hand, it turns out that many of the known interesting examples of non-separable Banach spaces (both in and outside ZFC) 
fall into CSP class.

\subsection*{General concepts and notations}

Real Banach spaces are typically denoted by $\X,\Y,\Z$. We denote by $\B_{\X}$ the closed unit ball of $\X$ and by $\S_{\X}$
the unit sphere of $\X$. Unless otherwise stated, we will apply cardinal arithmetic notations.

See \cite{En}, \cite{HHZ} and \cite{Ne} for the standard notions in set-theory, Banach spaces and topology, respectively.
We refer to Zizler's survey \cite{Z} on the non-separable Banach spaces for most of the definitions and results used here.

Recall that $\mathcal{F}\subset \X^{\ast}$ is a separating subset if and only if for each $x\in\X$ there is $f\in\mathcal{F}$
such that $f(x)\neq 0$ if and only if $\overline{\span}^{\omega^{\ast}}(\mathcal{F})=\X^{\ast}$ (see e.g. \cite[p.55]{HHZ}).
Let $\X$ be a Banach space and let $\mathcal{F}=\{(x_{\alpha},x^{\ast}_{\alpha})\}_{\alpha}\subset \X\times \X^{\ast}$ 
be a biorthogonal system, i.e. $x_{\alpha}^{\ast}(x_{\beta})=\delta_{\alpha,\beta}$. If $\overline{\span}(x_{\alpha})=\X$ and 
$\overline{\span}^{\omega^{\ast}}(x^{\ast}_{\alpha})=\X^{\ast}$, then $\mathcal{F}$ is called a \emph{Marku\v{c}evi\v{c}-basis} or 
\emph{M-basis}. Recall that each separable Banach space admits an M-basis, see \cite[p.219]{HHZ}. An M-basis 
$\{(x_{\alpha},f_{\alpha})\}_{\alpha}\subset \X\times \X^{\ast}$ is called countably $\lambda$-norming if there is $\lambda\geq 1$ 
such that for each $x\in \X$ there is $f\in\X^{\ast}$ satisfying $\lambda||x||\leq f(x)$ and 
$|\{\alpha:\ f(x_{\alpha})\neq 0\}|\leq \omega$. If $\X$ admits a countably norming M-basis, then $\X$ is called \emph{Plichko}
(reformulation according to \cite{Valdivia}). It is said that $\X$ has the \emph{density property} (DENS) if 
$\omega^{\ast}\mathrm{-dens}(\X^{\ast})=\mathrm{dens}(\X)$.

A compact Hausdorff space $K$ is called a \emph{Corson compact} if it can be embedded in a $\Sigma$-product of real lines, 
and a Banach space $\X$ is Weakly Lindel\"{o}f Determined (WLD) if $(\B_{\X^{\ast}},\omega^{\ast})$ is a Corson compact.
To mention a more general notion, recall that a compact Hausdorff space $K$ is a \emph{Valdivia compact} if there is 
an embedding $h\colon K\rightarrow \R^{\Gamma}$ (where $\R^{\Gamma}$ is equipped with the product topology) such that 
$h(K)\cap A$ is dense in $h(K)$, where $A\subset \R^{\Gamma}$ is the corresponding $\Sigma$-product. We will call a topological space 
\emph{countably perfect} if each point is in the countable closure of the rest of the space.

The following folklore facts will be applied frequently: 
Suppose that $(T,\tau)$ is a countably tight topological space, $\lambda$ is an ordinal 
with $\mathrm{cf}(\lambda)>\omega$ and $\{E_{\alpha}\}_{\alpha<\lambda}$ is a family of closed subsets of $T$
such that $E_{\alpha}\subset E_{\beta}$ for $\alpha<\beta$. 
Then $\overline{\bigcup_{\alpha<\lambda}E_{\alpha}}^{\tau}=\bigcup_{\alpha<\lambda}E_{\alpha}$.
A subset $\mathcal{F}\subset \X^{\ast}$ separates $X$ if and only if $\overline{\span}^{\omega^{\ast}}(\mathcal{F})=\X^{\ast}$
(see e.g. \cite[p.55]{HHZ}).

\begin{fact}\label{quotientfact}
Let $q\colon \X\rightarrow \X\quotient\Y$ be the quotient map $q\colon x\mapsto x+\Y$. Then 
$\overline{q(A)}=\{z+\Y|z\in\overline{A+\Y}\}$ for any subset $A\subset\X$.
\begin{proof}
The condition that $x\in \overline{A+\Y}$ is equivalent to
\[\inf_{y\in\Y,\ a\in A}||x-(a+y)||=\mathrm{dist}_{\X\quotient\Y}(q(x),q(A))=0.\]
\end{proof}
\end{fact} 
\newpage
\section{General properties}

It is easy to see that the following formulations of CSP are equivalent:
\begin{equation}\label{eq: basicchar}
\
\end{equation}
\begin{itemize}
\item{$\X$ has CSP: if $\mathcal{F}\subset\X^{\ast}$ separates $\X$, then there is a countable $\mathcal{F}_{0}\subset\mathcal{F}$, 
which separates $\X$.}

\item{For each $\omega^{\ast}$-dense linear subspace $V\subset\X^{\ast}$ there exists a countable subset $\mathcal{F}_{0}\subset V$
such that $\overline{\span}^{\omega^{\ast}}(\mathcal{F}_{0})=\X^{\ast}$.}
\item{Each family of closed subspaces with trivial intersection has a countable subfamily with trivial intersection.}

\item{There does not exist an uncountable separating family $\mathcal{F}\subset\X^{\ast}$ such that
each separating subfamily $\mathcal{F}_{0}\subset\mathcal{F}$ has the same cardinality as $\mathcal{F}$.}
\end{itemize}

Note that in particular $\omega^{\ast}\mathrm{-dens}(\X^{\ast})=\omega$ if $\X$ has CSP.

The spaces $c_{0}(\Gamma),\ \ell^{p}(\Gamma),\ 1\leq p\leq \infty,\ |\Gamma|\geq \omega_{1},$ provide examples of spaces 
without CSP, since $\{e_{\gamma}^{\ast}\}_{\gamma\in \Gamma}$ is an uncountable minimal separating family.

The following fact appeared already in \cite{BM} with a different formulation. 
\begin{proposition}\label{propo1}
Separable Banach spaces have CSP.
\end{proposition}
This follows also immediately from the following observation. 
\begin{proposition}
Let $\X$ be a Banach space such that $\S_{\X}$ is Lindel\"of in the relative weak topology. Then $X$ has CSP. 
\begin{proof}
Let $\{f_{\gamma}\}_{\gamma\in \Gamma}\subset \X^{\ast}$ be a separating family of functionals. 
Then $\bigcup_{\gamma\in\Gamma}f_{\gamma}^{-1}(\R\setminus\{0\})=\X\setminus \{0\}$. In particular,
the subsets $U_{\gamma}=f_{\gamma}^{-1}(\R\setminus\{0\})\cap \S_{\X},\ \gamma\in\Gamma,$ define a $\omega$-open cover for $\S_{\X}$.
According to the Lindel\"{o}f property of $(\S_{\X},\omega)$ there exists a countable subcover 
$\{U_{\gamma_{n}}\}_{n<\omega}\subset \{U_{\gamma}\}_{\gamma\in\Gamma}$. We conclude that $\{f_{\gamma_{n}}\}_{n<\omega}$ 
is a separating family of $X$.  
\end{proof}
\end{proposition}
We will subsequently give some examples of non-separable CSP spaces. It is straightforward to verify that CSP is preserved 
under isomorphisms.

\begin{proposition}\label{subsp}
Let $\X$ be a Banach space with CSP and let $\Y\subset \X$ be a closed subspace. Then $\Y$ has CSP.
\begin{proof}
Assume to the contrary that $\Y$ fails CSP. Suppose that $\mathcal{F}\subset \Y^{\ast}$ is an uncountable minimal 
separating set for $\Y$. Let $\widetilde{\mathcal{F}}\subset \X^{\ast}$ be a set of functionals obtained from 
$\mathcal{F}$ by Hahn-Banach extension. Then $\Y^{\bot}\cup \widetilde{\mathcal{F}}$ is an uncountable minimal 
separating subfamily for $\X$; a contradiction.
\end{proof} 
\end{proposition}

\begin{example}
\textit{The space $\ell^{\infty}$ does not have CSP.} 
\end{example}
The space $\ell^{\infty}$ contains an isometric copy of $\ell^{1}(2^{\omega})$ (see \cite[p.86]{HHZ}). 
Clearly $\ell^{1}(2^{\omega})$ does not have CSP. Thus the claim follows by applying Proposition \ref{subsp} that CSP is 
inherited by the closed subspaces.  

The preceding example shows that there is a countable $\mathcal{F}_{0}\subset (\ell^{\infty})^{\ast}$
such that 
\[ \span(\{f|_{c_{0}}:\ f\in \mathcal{F}_{0}\})\ \mathrm{is\ dense\ in}\ (\ell^{1},\omega^{\ast}),\]
and 
\[\span(\mathcal{F}_{0})\ \mathrm{is\ not\ dense\ in}\ ((\ell^{\infty})^{\ast},\omega^{\ast}),\]
even though $\ell^{1}$ is $\omega^{\ast}$-dense in $(\ell^{\infty})^{\ast}$ by Goldstine's theorem.
 
\begin{example}\label{JLc}
The spaces $JL_{0},JL_{2}$ (see \cite{JL}) have CSP according to subsequent observations (see Proposition \ref{propo}), but 
$JL_{0}\quotient c_{0}=c_{0}(2^{\omega})$ and $JL_{2}\quotient c_{0}=l^{2}(2^{\omega})$ (see e.g. \cite[p. 1757]{Z}) clearly do not. 
We conclude that
\begin{enumerate}
\item[(i)]{CSP does not pass to quotients in general.}
\item[(ii)]{The dual of a CSP space may contain a $\omega^{\ast}$-closed subspace, which is not $\omega^{\ast}$-separable.}
\end{enumerate}
The spaces $JL_{0}$ and $JL_{2}$ are not Lindel\"{o}f in their weak topology (see e.g. \cite[p.1757,1764]{Z}).
\end{example}

If $\X$ and $\Y$ have CSP, does it follow that $\X\oplus\Y$ has CSP? If so, is CSP a three-space property, i.e.
does $\X$ have CSP whenever $X\quotient \Y$ and $\Y\subset \X$ have CSP?
We will give a partial answer to this problem in Theorem \ref{Thm3SP}. 

Even though CSP is not a sufficient condition for separability, it is still quite a strong condition 'towards separability'.
\begin{proposition}\label{predual}
Let $\X$ be a Banach space such that $\X^{\ast}$ has CSP. Then $\X$ is separable.
\begin{proof}
Clearly $\X\subset \X^{\ast\ast}$ embedded canonically separates $\X^{\ast}$. Then according to CSP one can find a sequence 
$(x_{n})_{n\in N}\subset \X$ separating $\X^{\ast}$. We claim that $\overline{\span}(x_{n})=\X$. Indeed, if this was not the case, 
then one could find by the Hahn-Banach theorem a non-zero functional $f\in \X^{\ast}$ vanishing on $\overline{\span}(x_{n})$. 
But this is not possible, since $(x_{n})\subset \X^{\ast\ast}$ separates $\X^{\ast}$. Thus $\X$ is separable.
\end{proof}
\end{proposition} 
 
As a brief remark we would like to mention a mild condition under which separability and CSP are equivalent.

Let us consider the case that $\X$ admits a system 
$(\{x_{\alpha}\}_{\alpha},\{f_{\beta}\}_{\beta})\in \mathcal{P}(\X)\times \mathcal{P}(\X^{\ast})$ satisfying the following conditions:
\begin{equation}\label{eq: SEMIDENS}
\left\{ \begin{array}{rl}
 & \overline{\span}(\{x_{\alpha}\}_{\alpha})=\X\\
 & \overline{\span}^{\omega^{\ast}}(\{f_{\beta}\}_{\beta})=\X^{\ast}\\
 & \mathrm{For\ each}\ \beta\ \mathrm{we\ have}\ |\{\alpha |f_{\beta}(x_{\alpha})\neq 0\}|\leq \omega^{\ast}\mathrm{-dens}(\X^{\ast}).\end{array} \right.
\end{equation}
For a Banach space $\X$ the $\mathrm{DENS}$ property or the existence of an M-basis clearly provide 
a system satisfying \eqref{eq: SEMIDENS}.

\begin{proposition}\label{remarkSEMI}
A Banach space $\X$ is separable if and only if it has CSP and admits a system satisfying \eqref{eq: SEMIDENS}. 
In particular each CSP space with an M-basis is separable.
\end{proposition}
\begin{proof}
Suppose that $\X$ is a CSP space, and the system $(\{x_{\alpha}\}_{\alpha},\{f_{\beta}\}_{\beta})$ satisfies \eqref{eq: SEMIDENS}.
Then there exists a countable separating subfamily $\mathcal{F}\subset \{f_{\beta}\}_{\beta}$, since $\X$ has CSP.
Thus 
\[|\{\alpha|x_{\alpha}\neq 0\}|=|\{\alpha|f(x_{\alpha})\neq 0\ \mathrm{for\ some}\ f\in \mathcal{F}\}|\leq |\mathcal{F}\times (\omega^{\ast}\mathrm{-dens}(\X^{\ast}))|=\omega.\] 
Hence $\X$ is separable. On the other hand, each separable Banach space $\X$ has an M-basis (see e.g. \cite[p.219]{HHZ}) and 
hence system \eqref{eq: SEMIDENS}. Recall that, by Proposition \ref{propo1}, separability of $\X$ implies CSP.
\end{proof}

\begin{example}\label{Plichko}
\textit{For a Valdivia compact $K$ the space $C(K)$ has CSP if and only if $C(K)$ is separable.
If a Banach space $\X$ is Plichko and has CSP, then $\X$ is separable.}
\end{example}
Indeed, if $K$ is Valdivia, then $C(K)$ is Plichko (see \cite{Kalenda}). Any Plichko space admits an M-basis,
so that Proposition \ref{remarkSEMI} yields that $\X$ is separable.

\begin{proposition}
Let $\X$ be a CSP space and $C\subset \X$ a weakly compact set. Then $C$ is separable.
\end{proposition}
\begin{proof}
Observe that $\Y=\overline{\span}(C)$ is a WCG subspace and hence admits an M-basis. According to Proposition \ref{subsp}
the space $\Y$ has CSP. Thus, by Proposition \ref{remarkSEMI} $\Y$ is separable and so is $C$.
\end{proof} 

\section{Topological point of view}

A topological space $T$ is called \emph{dense-separable} if each dense subset $A\subset T$ is separable (see \cite{LM} for discussion).

Recall the following concept due to Corson and Pol. A Banach space $\X$ has property $(C)$ if, 
for each family $\mathcal{A}$ of closed convex subsets of $\X$ having empty intersection, there exists a countable subfamily 
$\mathcal{A}_{0}\subset\mathcal{A}$ with empty intersection. 
Pol gave an important characterization of property $(C)$ in terms of a dual space formulation $(C^{\prime})$, 
which is a kind of convex version of countable tightness (see \cite[Thm.3.4]{Pol}). This condition appears in the result below:

\begin{proposition}\label{propo}
Let $\X$ be a Banach space. Consider the following conditions:
\begin{enumerate}
\item[(1)]{$(\X^{\ast},\omega^{\ast})$ is countably tight i.e. for each $a\in \overline{A}^{\omega^{\ast}}\subset \X^{\ast}$
there is a subset $(a_{n})_{n<\omega}\subset A$ such that $a\in\overline{\{a_{n}|n<\omega\}}^{\omega^{\ast}}$.}
\item[(2)]{$\X^{\ast}$ is dense-separable in the $\omega^{\ast}$-topology.}
\item[(3)]{$\X$ satisfies property $(C^{\prime})$: for each $A\subset \X^{\ast}$ and $f\in\overline{A}^{\omega^{\ast}}$ there is a countable 
$A_{0}\subset A$ such that $f\in \overline{\conv}^{\omega^{\ast}}(A_{0})$.}
\item[(4)]{For each $A\subset \X^{\ast}$ and $f\in\overline{A}^{\omega^{\ast}}$ there is countable $A_{0}\subset A$ such that $f\in \overline{\span}^{\omega^{\ast}}(A_{0})$.}
\item[(5)]{$\X$ has CSP i.e. for each $A\subset \X^{\ast}$ such that $\overline{\span}^{\omega^{\ast}}(A)=\X^{\ast}$
there is a countable $A_{0}\subset A$ such that $\overline{\span}^{\omega^{\ast}}(A_{0})=\X^{\ast}$.} 
\end{enumerate}
If $\X^{\ast}$ is $\omega^{\ast}$-separable, then $1\implies 2\implies 5$ and $1\implies 3\implies 4\implies 5$.
\begin{proof}
Let us first check that implication $(1)\implies (2)$ holds if $\X^{\ast}$ is $\omega^{\ast}$-separable. 
Let $\X^{\ast}$ be $\omega^{\ast}$-separable space satisfying $(1)$ and let $A\subset \X^{\ast}$ be a $\omega^{\ast}$-dense subset. 
Fix a $\omega^{\ast}$-dense subset $(x_{k})_{k<\omega}\subset \X^{\ast}$. By using the countable tightness of $(\X^{\ast},\omega^{\ast})$, 
we may pick $\{a_{k}^{(n)}|n,k<\omega\}\subset A$ such that $x_{k}\in \overline{\{a_{k}^{(n)}|n<\omega\}}^{\omega^{\ast}}$ for each 
$k<\omega$. Hence $\overline{\{a_{k}^{(n)}|n,k<\omega\}}^{\omega^{\ast}}=\X^{\ast}$, so that $(2)$ holds, as $A$ was arbitrary.  

Let us check the implication $(4)\implies (5)$ for $\omega^{\ast}$-separable $\X^{\ast}$. First recall \eqref{eq: basicchar}. 
Let $\X^{\ast}$ be $\omega^{\ast}$-separable space satisfying $(4)$ and let 
$(x_{n})_{n<\omega}\subset \X^{\ast}$ be a $\omega^{\ast}$-dense subset in $\X^{\ast}$. Consider a separating family 
$A\subset \X^{\ast}$. Thus $\span(A)$ is $\omega^{\ast}$-dense in $\X^{\ast}$. According to condition (4) we can find a countable set 
$C_{n}\subset \span(A)$ for each $n<\omega$ such that $x_{n}\in \overline{\span}^{\omega^{\ast}}(C_{n})$. 
Thus $\overline{\span}^{\omega^{\ast}}(\bigcup_{n} C_{n})=\X^{\ast}$.
Note that each of the sets $C_{n}$ is contained in the linear span of countably many vectors of $A$. Therefore
there is countable $A_{0}\subset A$ such that 
\[\overline{\span}^{\omega^{\ast}}(A_{0})=\overline{\span}^{\omega^{\ast}}(\bigcup_{n}C_{n})=\X^{\ast},\]
so that $\X$ satisfies $(5)$. Other implications are clear.
\end{proof}
\end{proposition}

We do not know if $(2)\implies (3)$ above or if some implications can be reversed.

For example the following spaces are non-separable and have CSP according to Proposition \ref{propo}: 
$JL_{0},\ JL_{2}$ and $C(K)$, where $K$ is the double-arrow space. Indeed, these spaces have property $(C)$ and a 
$\omega^{\ast}$-separable dual (see e.g. \cite[p.1757]{Z}, \cite[p.146]{Pol}).

Assuming CH Kunen constructed an interesting Hausdorff compact $K$, which is separable, scattered and non-metrizable 
(see \cite[p.1086-1099]{Ne} for discussion). Kunen's $C(K)$ space is non-separable and $(C(K)^{\ast},\omega^{\ast})$ is 
hereditarily separable (see \cite[p.476]{HSZ}), hence separable and countably tight.

\begin{example}
\textit{There is a bounded injective linear operator $T\colon JL_{0}\rightarrow c_{0}(2^{\omega})$ whose range is non-separable.}
\end{example}
In justifying this we will apply the fact that $JL_{0}$ contains $c_{0}$ and that $JL_{0}\quotient c_{0}=c_{0}(2^{\omega})$ (see \cite[p. 1757]{Z}).
Let $(e_{n})_{n<\omega}\subset c_{0}$ be the canonical unit vector basis and let $(e_{n}^{\ast})_{n<\omega}\subset \ell^{1}$
be the corresponding functionals. Let $(f_{n})_{n<\omega}\subset JL_{0}^{\ast}$ be the Hahn-Banach extension of 
$(e_{n})_{n<\omega}$. Define $S\colon JL_{0}\rightarrow c_{0}$ by $x\mapsto ((n+1)^{-1}f_{n}(x))_{n<\omega}$. Let 
$q\colon JL_{0}\rightarrow JL_{0}\quotient c_{0}$ be the canonical quotient mapping. Then 
$T\colon x\mapsto (S(x),q(x))$ defines the required map 
$JL_{0}\rightarrow c_{0}(\omega)\oplus c_{0}(2^{\omega})=c_{0}(2^{\omega})$.

\begin{problem}
Suppose that $\X$ admits a long unconditional basis and $\Y\subset\X$ is a closed subspace having CSP. Does it follow that $\Y$ is 
separable?
\end{problem}

Dual CSP spaces can be characterized as follows:
\begin{theorem}\label{dualchar}
Let $\X$ be a Banach space. Then the following are equivalent:
\begin{enumerate}
\item{$\X^{\ast}$ has CSP}
\item{$\X$ is separable and $\X^{\ast}$ has property (C)}
\item{$\X$ is separable and does not contain $\ell^{1}$ isomorphically.}
\end{enumerate}
\begin{proof}
The equivalence of the last two conditions is known (see \cite[Thm.4.2]{Z}). 

If $(2)$ holds, then $\X^{\ast\ast}$ is $\omega^{\ast}$-separable by Goldstine's theorem, so that Proposition \ref{propo} 
can be applied together with property $(C)$ to obtain that $\X^{\ast}$ has CSP.

On the other hand, if $\X^{\ast}$ has CSP, then by Proposition \ref{predual} we known that $\X$ must be separable. 
Now assume to the contrary that $\X$ contains an isomorphic copy of $\ell^{1}$. 
Then it is known that $\X^{\ast}$ contains a complemented isomorphic copy of $\ell^{\infty}$. 
Since $\ell^{\infty}$ does not have CSP it follows by Proposition \ref{subsp} 
that $\X^{\ast}$ fails CSP, a contradiction. Thus $\X$ does not contain $\ell^{1}$.
\end{proof}
\end{theorem}

For example the James Tree and the James function space (see \cite{LS}) are separable spaces not containing $\ell^{1}$, and whose duals
$JT^{\ast}$ and $JF^{\ast}$ are non-separable CSP spaces.

\subsection{$C(K)$ spaces}

\begin{proposition}\label{ast}
Let $L$ be a locally compact Hausdorff space such that $C_{0}(L)$ has CSP. 
Then $L$ is dense-separable and the interior of the derived set of $L$ is countably perfect.
\begin{proof}
Recall that locally compact Hausdorff spaces are completely regular.
If the interior of the derived set of $L$ is non-empty, then let $x\in \mathrm{int}(D(L))$.
If such $x$ exists, it is not an isolated point in $L$, and in any case we may fix a dense subset $\Gamma\subset L$ such that 
$x\notin \Gamma$. Consider the point evaluation maps
$f\stackrel{\delta_{k}}\mapsto f(k)$ in $C_{0}(L)^{\ast}$ for $k\in \Gamma$. 
Clearly these maps separate $C_{0}(L)$. Since $C_{0}(L)$ has CSP, there exists a countable subset $\Gamma_{0}\subset \Gamma$ 
such that the associated evaluation maps still separate $C_{0}(L)$. We claim that $L=\overline{\Gamma}_{0}$. Indeed, if there 
exists a point $y\in L\setminus \overline{\Gamma}_{0}$, then by the completely regularity of $L$ there exists
$f\in C_{0}(L)$ attaining value $1$ at $y$ but vanishing on $\overline{\Gamma}_{0}$. This contradicts the fact 
that the point evaluations associated to the points in $\Gamma_{0}$ separate $C_{0}(L)$. 
Thus $\overline{\Gamma}_{0}=L$ and $L$ is dense-separable as $\Gamma$ was arbitrary. 
Note that $x\in \overline{\Gamma_{0}\cap \mathrm{int}(D(L))}$, so that $\mathrm{int}(D(L))$ is countably perfect.
\end{proof}
\end{proposition}

We have not been able to construct a (dense-separable) compact $K$ such that $C(K)$ has CSP but fails property (C).
Observe that the \v{C}ech-Stone compactification $\beta\omega$ is dense-separable as $n\in\omega$ are isolated in $\beta\omega$. 
Since $C(\beta\omega)=\ell^{\infty}$, we conclude that dense-separability of $K$ is not sufficient for $C(K)$ to have CSP.

The following result produces plenty of examples of non-separable CSP spaces.
\begin{theorem}\label{sca}
Let $K$ be a scattered and countably tight compact. Then $C(K)$ has CSP if and only if $K$ is separable.
\begin{proof}
The only if part follows from Proposition \ref{ast}. Towards the other implication, recall that 
a scattered compact $K$ is countably tight if and only if $C(K)$ has property (C) (see \cite[Cor.4.1]{Pol}).
If $K$ is separable, then a standard argument using the point evaluations gives that $C(K)^{\ast}$ is $\omega^{\ast}$-separable, 
so that Proposition \ref{propo} can be applied.
\end{proof}
\end{theorem}
For example the one-point compactified rational sequence topology (see \cite[p.87]{CET}) is scattered, countably tight,
separable and non-metrizable. 

Analogous to the open question about preservation of CSP in finite sums is the open question about the preservation of 
dense-separability in finite products of topological spaces. The following fact, however, is easy to verify.
\begin{remark}
Let $A$ and $B$ be topological spaces. If $A$ is dense-separable and $B$ has a countable $\pi$-base, then
$A\times B$ is dense-separable.
\end{remark} 

\subsection{Auxiliary results: Intersections of distended subspaces}

In addition to property $(C)$ we will treat another type of condition concerning intersections of convex sets, 
which seems to be closely related to $(C)$. Throughout this subsection let $\X$ be a Banach space and $\kappa$ 
an uncountable regular ordinal. We denote here by $\{Z_{\sigma}\}_{\sigma<\kappa}$ a family of closed subspaces of $\X$ such that 
$Z_{\alpha}\supsetneq Z_{\beta}$ for $\alpha<\beta<\kappa$ and $\bigcap_{\sigma<\kappa}Z_{\sigma}=\{0\}$ and
we impose the existence of such family for $\X$. This actually excludes $\X$ from CSP class, as it will turn out in the next section.

Some subsequent results here depend on the following question:\\ 
\textit{Given $\X$ and $\{Z_{\sigma}\}_{\sigma<\kappa}$ as above, is $\bigcap_{\sigma<\kappa}(\B_{\X}+Z_{\sigma})$ bounded?}\\
At first glance the answer may appear to be positive for all Banach spaces $\X$. For example, it is easy to see that if $\X$ is reflexive, 
then $\bigcap_{\sigma<\kappa}(\B_{\X}+Z_{\sigma})=\B_{\X}$. However, next we present an example of a space for which
the answer to the above question is negative. Define a function $|||\cdot|||\colon \ell^{\infty}(\omega_{1})\rightarrow [0,\infty]$
by 
\[|||(x_{\alpha})_{\alpha<\omega_{1}}|||=||(x_{\alpha})_{\alpha<\omega_{1}}||_{\ell^{\infty}(\omega_{1})}+\sum_{n<\omega}n\limsup_{i\rightarrow\omega_{1}}|x_{\omega i+n}|\quad (\mathrm{ordinal\ arithmetic}).\]
Then $(\X,|||\cdot|||)$, where $\X=\{x\in\ell^{\infty}(\omega_{1}):\ |||x|||<\infty\}$, is clearly a Banach space.
Let $E_{\sigma}=\{ (x_{\alpha})_{\alpha<\omega_{1}} \in \X :\ x_{\alpha}=0\ \mathrm{for}\ \alpha<\sigma\}$ for $\sigma<\omega_{1}$.
We denote by $\1_{A}\colon [0,\omega_{1})\rightarrow \{0,1\}$ the characteristic function of a given subset $A\subset \omega_{1}$.
Note that $\1_{[0,\sigma]}\in \S_{(\X,|||\cdot|||)}$ for all $\sigma<\omega_{1}$.
Hence $\1_{\{\omega i +n|i<\omega_{1}\}}\in\bigcap_{\sigma<\omega_{1}}(\B_{(\X,|||\cdot|||)}+E_{\sigma})$ for all $n<\omega$.
Consequently $\bigcap_{\sigma<\omega_{1}}(\B_{(\X,|||\cdot|||)}+E_{\sigma})$ is unbounded.

For convenience we denote by ($\mathcal{B}$) the class of Banach spaces $\X$ satisfying that $\bigcap_{\sigma<\kappa}\B_{\X}+Z_{\sigma}$
is bounded for any $\{Z_{\sigma}\}_{\sigma<\kappa}$ such as above (or trivially if such $\{Z_{\sigma}\}_{\sigma<\kappa}$ 
does not exist at all).
\begin{proposition}
For any $\X$ we have that $\bigcap_{\epsilon>0}\bigcap_{\sigma<\kappa}(\epsilon\B_{\X}+Z_{\sigma})=\{0\}$.
\end{proposition}
\begin{proof}
Observe that $\bigcap_{\sigma<\kappa}(\epsilon\B_{\X}+Z_{\sigma})$ is a symmetric convex set, which 
contains $0$.

First we wish to check that $\bigcap_{\epsilon>0}\bigcap_{\sigma<\kappa}(\epsilon\B_{\X}+Z_{\sigma})$ does not contain
any non-trivial linear subspace $L$. Indeed, if $\bigcap_{\sigma<\kappa}(\B_{\X}+Z_{\sigma})$ contains $L=[x]$ 
for some $x\in\S_{\X}$, then $\dist(kx,Z_{\sigma})\leq 1$ for all $k<\omega$ and $\sigma<\kappa$.
Observe that 
\begin{equation}\label{eq: epsilonB}
\epsilon\B_{\X}+Z_{\sigma}=\epsilon\B_{\X}+\epsilon Z_{\sigma}=\epsilon(\B_{\X}+Z_{\sigma})
\end{equation}
for all $\epsilon>0$ and $\sigma<\kappa$. Thus, by putting $\epsilon=\frac{1}{k}$ we obtain that 
$\dist(x,Z_{\sigma})=\frac{\dist(kx,Z_{\sigma})}{k}<\frac{1}{k}$ for all $k<\omega$ and $\sigma<\kappa$.
Since $Z_{\sigma}$ are closed subspaces, we get that $x\in Z_{\sigma}$ for all $\sigma<\kappa$. 
This contradicts the fact that $\bigcap_{\sigma<\kappa}Z_{\sigma}=\{0\}\not\ni x$ and hence $\X$ does not contain 
any non-trivial linear subspace $L$.

Now, let $L$ be a $1$-dimensional subspace and write $l=L\cap \bigcap_{\sigma<\kappa}(\B_{\X}+Z_{\sigma})$.
Since $\bigcap_{\sigma<\kappa}(\B_{\X}+Z_{\sigma})$ is convex, symmetric and does not contain $L$, we obtain
that $l$ is bounded. Hence $L\cap \bigcap_{\epsilon>0}\epsilon\bigcap_{\sigma<\kappa}(\B_{\X}+Z_{\sigma})=\{0\}$.
Since $L$ was arbitrary, we obtain that $\bigcap_{\epsilon>0}\bigcap_{\sigma<\kappa}(\epsilon\B_{\X}+Z_{\sigma})=\{0\}$.
\end{proof}
Note that in the above proof we did not need that $\kappa$ is uncountable or regular. 
Next we will give results towards applications of subsequent Lemma \ref{cap}, which is our main technical machinery. 

\begin{proposition}
Suppose that the dual of $\X$ satisfies the condition $(4)$ of Proposition \ref{propo}. Then $\X$ is a member of ($\mathcal{B}$).
\end{proposition}
\begin{proof}
Let $\{Z_{\sigma}\}_{\sigma<\kappa}$ be a strictly nested sequence of closed subspaces of $\X$ such that 
$\bigcap_{\sigma<\kappa}Z_{\sigma}=\{0\}$. Observe that 
$\overline{\span}^{\omega^{\ast}}(\bigcup_{\sigma<\kappa}Z_{\sigma}^{\bot})=
\overline{\bigcup_{\sigma<\kappa}Z_{\sigma}^{\bot}}^{\omega^{\ast}}=\X^{\ast}$.
Fix $f\in \X^{\ast}$. Since $\kappa$ is regular, condition $(4)$ of Proposition \ref{propo} yields that 
there is $\sigma_{0}<\kappa$ such that $f\in Z_{\sigma}^{\bot}$ whenever $\sigma_{0}\leq \sigma<\kappa$.
Thus $f(\bigcap_{\sigma<\kappa}\B_{\X}+Z_{\sigma})\subset [-||f||,||f||]$. Note that $\bigcap_{\sigma<\kappa}\B_{\X}+Z_{\sigma}$
is weakly bounded as $f$ was arbitrary. Hence $\bigcap_{\sigma<\kappa}\B_{\X}+Z_{\sigma}$ is bounded according to the 
uniform boundedness principle.
\end{proof}
  
\begin{proposition}
Suppose that the dual of $\X$ satisfies the condition $(4)$ of Proposition \ref{propo}, that $\{Y_{\sigma}\}_{\sigma<\kappa}$ is a nested 
family of affine subspaces, where $\kappa$ is an uncountable regular cardinal, and that $\bigcap_{\sigma<\kappa}Y_{\sigma}=\{y\}$.
If $y_{\sigma}\in Y_{\sigma}$ for $\sigma<\kappa$, then $y_{\sigma}\rightarrow y$ weakly as $\sigma\rightarrow \kappa$.
\end{proposition}
\begin{proof}

Suppose that $u+Z_{1}=Y_{1}\supset Y_{2}=u+v+Z_{2}$ are closed affine subspaces, 
where $u,v\in \X$ and $Z_{1},Z_{2}\subset \X$ are closed subspaces. 
Then $v\in Z_{\sigma_{1}}$ and $Z_{\sigma_{2}}\subset Z_{\sigma_{1}}$.
Indeed, if $Z_{1}\not\supset Z_{2}$, then there is no $u+v\in\X$ such that $u+Z_{1}\supset u+v+Z_{2}$.
On the other hand, if $v\notin Z_{1}$, then $v+Z_{2}\not\subset Z_{1}$, so that 
$u+v+Z_{2}\not\subset u+\Z_{1}$. 

Let $\{Z_{\sigma}\}_{\sigma<\kappa}$ be the nested family of closed subspaces of $\X$ corresponding to 
$\{Y_{\sigma}\}_{\sigma<\kappa}$. It follows that $y_{\sigma}\in y+ Z_{\sigma}$ for $\sigma<\kappa$ and 
$\bigcap_{\sigma<\kappa}Z_{\sigma}=\{0\}$. 

Fix $f\in\X^{\ast}$. Observe that $\overline{\span}^{\omega^{\ast}}(\bigcup_{\sigma<\kappa}Z_{\sigma}^{\bot})=\X^{\ast}$ and 
$Z_{\sigma}^{\bot}$ is $\omega^{\ast}$-closed for $\sigma<\kappa$. By using condition (4) of
Proposition \ref{propo} and the fact that $\cf(\kappa)>\omega$, we obtain that there is $\sigma_{0}<\kappa$
such that $f\in Z_{\sigma_{0}}^{\bot}$. Then $f(y_{\sigma}-y)=0$ for $\sigma\geq \sigma_{0}$, so that we have the claim
as $f$ was arbitrary.
\end{proof} 
 
\begin{theorem}
Suppose that $\X$ satisfies the following condition: For each nested family of closed affine subspaces 
$\{Y_{\sigma}\}_{\sigma<\kappa}$ it holds that $\bigcap_{\sigma<\kappa}Y_{\sigma}\neq\emptyset$. 
Then $\X$ is a member of ($\mathcal{B}$).  
\end{theorem}
Before giving the proof, we note that $\sup_{\sigma<\kappa}\dist(0,Y_{\sigma})<\infty$, since $\kappa$ is 
uncountable and regular.
\begin{proof}
Let $\kappa$ and $\{Z_{\sigma}\}_{\sigma<\kappa}$ be as in the beginning of this subsection. 
Since $\bigcap_{\sigma<\kappa}Z_{\sigma}=\{0\}$, we 
may define a norm $|||\cdot|||\colon \X\rightarrow \R$ by $|||x|||=\sup_{\sigma<\kappa}\dist(x,Z_{\sigma})$. 
Observe that $|||x|||\leq ||x||$ for $x\in\X$. 

We aim to show that the norms $||\cdot||$ and $|||\cdot|||$ are equivalent. After this has been established, 
it follows that there is $C>0$ such that $||x||\leq C$ whenever $\dist(x,Z_{\sigma})\leq 1$ for $\sigma<\kappa$.
Actually, it suffices to check that $|||\cdot|||$ is complete. Indeed, in such case the Banach open mapping principle
yields that $\I\colon (\X,||\cdot||)\rightarrow (\X,|||\cdot|||)$ is an isomorphism.

Let $(x_{n})_{n<\omega}\subset \X$ be a $|||\cdot|||$-Cauchy sequence, where $x_{0}=0$. Denote
$\widehat{x}_{n}^{\sigma}=x_{n}+Z_{\sigma}\in \X\quotient Z_{\sigma}$ for $n<\omega,\ \sigma<\kappa$.
Observe that 
\[|||x_{i}-x_{j}|||=\sup_{\sigma<\kappa}||\widehat{x}^{\sigma}_{i}-\widehat{x}^{\sigma}_{j}||_{\X\quotient Z_{\sigma}}\quad 
\mathrm{for}\ i,j\in\omega.\]
Since $\cf(\kappa)>\omega$, there is for $(i,j)\in \omega\times\omega$ an ordinal $\alpha_{i,j}<\kappa$
such that 
\[|||x_{i}-x_{j}|||=||\widehat{x}^{\alpha_{i,j}}_{i}-\widehat{x}^{\alpha_{i,j}}_{j}||_{\X\quotient Z_{\alpha_{i,j}}}.\]
Observe that $\alpha\stackrel{\cdot}{=}\sup_{i,j\in\omega}\alpha_{i,j}<\kappa$. Next we regard $\sigma\in [\alpha,\kappa)$.

Thus the sequence $(\widehat{x}^{\sigma}_{n})_{n<\omega}\subset \X\quotient Z_{\sigma}$
is Cauchy for all $\sigma\in [\alpha,\kappa)$. Since $\X\quotient Z_{\sigma}$ is a Banach space, we obtain that there is
$y^{\sigma}\in \X\quotient Z_{\sigma}$ such that $\widehat{x}^{\sigma}_{n}\rightarrow y^{\sigma}$ as $n\rightarrow \infty$
in $\X\quotient Z_{\sigma}$. Moreover, by the selection of $\alpha$ we get 
$\lim_{n\rightarrow\infty}|||x_{n}|||=\lim_{n\rightarrow\infty}||\widehat{0}^{\sigma}-\widehat{x}^{\sigma}_{n}||_{\X\quotient Z_{\sigma}}$, since $x_{0}=0$. 

We may regard $y^{\sigma}\stackrel{\cdot}{=}Y_{\sigma}\stackrel{\cdot}{=}v_{\sigma}+Z_{\sigma}\subset \X$ as affine subspaces, 
where $v_{\sigma}\in \X$ for $\sigma\in [\alpha,\kappa)$. Similar interpretation for $\widehat{x}_{n}^{\sigma}\in\X$ 
yields the following: 
Since $\widehat{x}_{n}^{\sigma_{1}}=\widehat{x}_{n}^{\sigma_{2}}\quotient Z_{\sigma_{1}}$
for $n<\omega,\ \alpha\leq\sigma_{1}\leq \sigma_{2}<\kappa,$ we obtain that $y^{\sigma_{1}}=y^{\sigma_{2}}\quotient Z_{\sigma_{1}}$
by the continuity of the quotient map $\X\quotient Z_{\sigma_{2}}\rightarrow \X\quotient Z_{\sigma_{1}}$.
Thus $Y_{\sigma_{1}}=Y_{\sigma_{2}}+Z_{\sigma_{1}}$, where $\alpha\leq \sigma_{1}\leq \sigma_{2}<\kappa$. Note that 
$\dist(0,Y_{\sigma})=\lim_{n\rightarrow\infty}|||x_{n}|||$ for $\sigma\in [\alpha,\kappa)$ by the selection of $\alpha$. 
According to the assumptions 
$\bigcap_{\sigma<\kappa}Y_{\sigma}\neq\emptyset$ and let us choose $x\in \bigcap_{\sigma<\kappa}Y_{\sigma}$ 
(even though it turns out promptly that this set is a singleton). Observe that 
$y^{\sigma}=\widehat{x}^{\sigma}\stackrel{\cdot}{=}x+Z_{\sigma}$ for $\sigma\in [\alpha,\kappa)$. 
Thus $\sup_{\sigma\in [\alpha,\kappa)}||\widehat{x}^{\sigma}_{n}-\widehat{x}^{\sigma}||_{\X\quotient Z_{\sigma}}\rightarrow 0$ 
as $n\rightarrow\infty$. Hence $x_{n}\stackrel{|||\cdot|||}{\longrightarrow} x$ as $n\rightarrow \infty$, which
completes the proof.
\end{proof}
 
\begin{lemma}\label{cap}
Let $\X$ be a Banach space, $\Y\subset \X$ a closed subspace and $\kappa$ an uncountable regular cardinal. Let $Z_{\sigma}\subset \X$ 
be closed subspaces, for $\sigma<\kappa$, which satisfy $Z_{\alpha}\subsetneq Z_{\beta}$ for $\beta<\alpha<\kappa$
and $\bigcap_{\sigma<\kappa}Z_{\sigma}=\{0\}$. Then the following facts hold:
\begin{enumerate}
\item[(i)]{If $\dens(\Y)<\kappa$ and $\X\in (\mathcal{B})$, 
then $\Y=\bigcap_{\sigma<\kappa}\overline{\Y+Z_{\sigma}}$.}
\item[(ii)]{If $\dens(\Y)<\kappa$, then there exists $\theta<\kappa$ such that $Z_{\theta}\cap \Y=\{0\}$.}
\end{enumerate}
\end{lemma}
\begin{proof}
Let us treat the claim $(i)$. Let $x\in \bigcap_{\sigma<\kappa}\overline{\span}(\Y\cup Z_{\sigma})$ and $\epsilon>0$.
Thus there is a cofinal sequence $\{\sigma_{\alpha}\}_{\alpha<\kappa}\subset \kappa$ and families 
$\{y_{\alpha}\}_{\alpha<\kappa}\subset \Y$ and $\{z_{\alpha}\}_{\alpha<\kappa}\subset \X$
such that $z_{\alpha}\in Z_{\sigma_{\alpha}}$ and $||x-(y_{\alpha}+z_{\alpha})||<\frac{\epsilon}{2}$
for each $\alpha<\kappa$.

Note that we have $\bigcap_{\beta<\kappa}\overline{\{y_{\alpha}|\beta<\alpha<\kappa\}}\neq\emptyset$, because
$\overline{\{y_{\alpha}|\beta<\alpha<\kappa\}}_{\beta<\kappa}$ is a decreasing sequence of closed sets in $\Y$,
where $\kappa$ is regular, and the Lindel\"{o}f number of $\Y$ is less than $\kappa$.

Let $y^{(\epsilon)}\in \bigcap_{\beta<\kappa}\overline{\{y_{\alpha}|\beta<\alpha<\kappa\}}$ depending on the choice of $\epsilon$. 
Hence we may pick a cofinal sequence $\{\alpha_{\delta}\}_{\delta<\kappa}\subset \kappa$ such that 
$||y^{(\epsilon)}-y_{\alpha_{\delta}}||<\frac{\epsilon}{2}$ for each $\delta<\kappa$. 
This means that
\[||x-(y^{(\epsilon)}+z_{\alpha_{\delta}})||\leq ||y^{(\epsilon)}-y_{\alpha_{\delta}}||+||x-(y_{\alpha_{\delta}}+z_{\alpha_{\delta}})||<\epsilon\]
for $\delta<\kappa$. We get that $x-y^{(\epsilon)}\in z_{\alpha_{\delta}}+\epsilon\B_{\X}$ and in particular
\begin{equation}\label{eq: xyeps}
x-y^{(\epsilon)}\in Z_{\sigma_{\alpha_{\delta}}}+\epsilon\B_{\X}\ \mathrm{for}\ \delta<\kappa.
\end{equation}
Since $\{Z_{\sigma}\}_{\sigma<\kappa}$ decreases to $\{0\}$ and $\sup_{\delta<\kappa}\alpha_{\delta}=\kappa$, 
we obtain that $\bigcap_{\delta<\kappa}Z_{\sigma_{\alpha_{\delta}}}=\{0\}$.

Since $\X$ belongs to ($\mathcal{B}$), it follows by \eqref{eq: epsilonB} that  
\[\lim_{\epsilon\rightarrow 0^{+}}\diam(\bigcap_{\delta<\kappa}(\epsilon\B_{\X}+Z_{\sigma_{\alpha_{\delta}}}))=0.\]
Since $\epsilon>0$ was arbitrary in \eqref{eq: xyeps}, we conclude  that $\dist(x,\Y)=0$, that is $x\in \Y$ as $\Y$ is closed.
This completes the proof of claim $(i)$. 

Let us check claim $(ii)$. Since $\mathrm{dens}(\Y)<\kappa$, the Lindel\"{o}f number of $\Y\setminus \{0\}$ is less than $\kappa$.
It follows by the regularity of $\kappa$, that there cannot exist a decreasing sequence 
$\{(Z_{\alpha}\cap \Y)\setminus \{0\}\}_{\alpha<\kappa}$ of non-empty closed sets in $\Y\setminus \{0\}$. 
\end{proof}

We note that in the above lemma the additional assumptions in (i) cannot be removed. Indeed, consider closed
subspaces 
\[E_{\alpha}=\{(x_{i})\in\ell^{\infty}(\omega_{1})|x_{i}=0,\ \mathrm{for}\ i\leq \alpha\},\quad \alpha<\omega_{1}\] 
of $\ell^{\infty}(\omega_{1})$. Then $\bigcap_{\alpha<\omega_{1}}\ell^{\infty}_{c}(\omega_{1})+E_{\alpha}=\ell^{\infty}(\omega_{1})$.
Note that $\mathrm{dens}(\ell^{\infty}_{c}(\omega_{1}))=2^{\omega}$ and $\mathrm{dens}(\ell^{\infty}(\omega_{1}))=2^{\omega_{1}}$.

\section{Combinatorial approach to CSP spaces}

Recall that if $\X$ has CSP, then $\omega^{\ast}\mathrm{-dens}(\X^{\ast})=\omega$. This in turn implies
$\mathrm{dens}(\X)\leq |\X|\leq 2^{\omega}$. Recall that $2^{\omega}<\aleph_{\omega}$ is consistent with ZFC.
For this reason some of the results here, which involve assumption about the density of the space, can actually be thought of 
as consistency results.  

For a given Banach space $\X$ we define the \emph{cofinality range} of $\X$, $\mathrm{cr}(\X)$ for short, as 
the set of all infinite regular cardinals $\kappa$ satisfying that there exists a family $\{E_{\sigma}\}_{\sigma<\kappa}$ 
of closed subspaces of $\X$ such that $\bigcap_{\sigma<\kappa}E_{\sigma}=\{0\}$ and $E_{\alpha}\subsetneq E_{\beta}$ whenever 
$\alpha<\beta<\kappa$. Observe that each $\kappa\in \mathrm{cr}(\X)$ satisfies $\kappa\leq\mathrm{dens}(\X)$.

\begin{theorem}\label{subspaceseq}
Let $\X$ be a Banach space such that $\mathrm{dens}(\X)<\aleph_{\omega}$. Then $\X$ has CSP if and only if 
$\mathrm{cr}(\X)=\{\omega\}$.
\begin{proof}
Let us first check the easier 'only if' part. Suppose that for $\{E_{\sigma}\}_{\sigma<\kappa}$, as above, there does not exist 
$(\sigma_{n})_{n<\omega}$ such that $\bigcap_{n<\omega}E_{\sigma_{n}}=\{0\}$, or equivalently $\sup_{n<\omega}\sigma_{n}=\kappa$.
Note that $\bigcup_{\sigma<\kappa}E_{\sigma}^{\bot}\subset \X^{\ast}$ separates $\X$ by the Hahn-Banach theorem. 
Clearly this set has no countable separating subset, so that $\X$ fails CSP.  

To check the 'if' part assume that $\X$ satisfies $\mathrm{dens}(\X)<\aleph_{\omega}$ and $\X$ fails CSP. 
We aim to show that in such case $\mathrm{cr}(\X)\neq \{\omega\}$. 
According to \eqref{eq: basicchar} there is an uncountable separating family $\mathcal{F}\subset \X^{\ast}$, 
whose all separating subfamilies have the same cardinality, say $\kappa\geq \omega_{1}$. 
Write $\mathcal{F}=\{f_{\alpha}\}_{\alpha<\kappa}$.  
Put $F_{\sigma}=\bigcap_{\alpha<\sigma}\ker f_{\alpha}$ for all $\sigma<\kappa$. Clearly this gives a 
(not necessarily strictly) nested family of closed subspaces. Note that $\bigcap_{\sigma<\kappa}F_{\sigma}=\{0\}.$

Let $\phi(0)=0$ and we define recursively
\[\phi(\alpha)=\min\{\beta<\kappa:\ \bigcap_{\gamma<\alpha}\ker(f_{\phi(\gamma)})\not\subset \ker(f_{\beta})\}.\]
Thus, putting $E_{0}=\X$, $E_{\sigma}=\bigcap_{\alpha<\sigma}\ker(f_{\phi(\alpha)})$ for $0<\sigma<\kappa$ 
gives a strictly nested family and $\bigcap_{\sigma<\lambda}E_{\sigma}=\{0\}$ for some ordinal $\lambda\leq \kappa$.
Note that by the construction of $\phi$ it holds that $\{f_{\phi(\alpha)}\}_{\alpha<\lambda}$ is a separating family.

Since $\{f_{\alpha}\}_{\alpha<\kappa}$ does not have a separating subfamily of cardinality less than $\kappa$,
we conclude that $\lambda=\kappa$ above. As $\mathrm{dens}(\X)<\aleph_{\omega}$ we obtain that 
$\lambda=\kappa<\aleph_{\omega}$ is a regular cardinal and $\{E_{\alpha}\}_{\alpha<\lambda}$ is the required family witnessing that 
$\omega_{1}\leq\lambda=\kappa\in \mathrm{cr}(\X)$. 
\end{proof}
\end{theorem}

The above proof yields immediately that CSP$\implies(\mathcal{B})$. 
Note that there are CSP spaces with an uncountable biorthogonal system. 
For example, it suffices to consider a non-separable dual space with CSP, such as $JT^{\ast}$, since it is known that non-separable 
dual spaces have uncountable biorthogonal systems (see \cite[Cor.4]{St}).

One could ask if the previous result remains valid if one replaces 'cofinality' by 'cardinality' in the definition of $\mathrm{cr}(\X)$. 
This is not the case as the following example shows; the assumption about the regularity of $\kappa$ in the definition of 
$\mathrm{cr}(\X)$ is indeed essential:

\begin{example}
If $\X$ is a CSP space with a biorthogonal system $\{(x_{\alpha},x_{\alpha}^{\ast})\}_{\alpha<\omega_{1}}$, 
then by using the proof of the previous theorem we see that
\begin{enumerate}
\item{$F\stackrel{\cdot}{=}\bigcap_{\beta<\omega_{1}}\overline{\span}(\{x_{\alpha}|\beta<\alpha<\omega_{1}\})\neq \{0\}$}
\item{$\X\quotient F$ is not a CSP space.}
\item{There is a strictly nested family $\{E_{\sigma}\}_{\sigma<\kappa}$ of closed subspaces of $\X$ such that
$\bigcap_{\sigma<\kappa}E_{\sigma}=\{0\}$, where $\kappa>\omega_{1}$.}
\end{enumerate}
\end{example}
The family of subspaces in the last condition can be obtained as follows: Let $\{F_{\theta}\}_{\theta<\lambda}$
be a strictly nested family of closed subspaces of $F$, where $\lambda$ is an ordinal and $\bigcap_{\theta<\lambda}F_{\theta}=\{0\}$.
It suffices to put $E_{\sigma}=\overline{\span}(\{x_{\alpha}|\sigma\leq \alpha<\omega_{1}\})$ for $\sigma<\omega_{1}$,
$E_{\omega_{1}}=F$, and $E_{\omega_{1}+\theta}=F_{\theta+1}$ for $\theta<\lambda$.  

We attempt to bring CSP closer to a three space property in the next result.

\begin{theorem}\label{Thm3SP}
Let $\X$ be a Banach space satisfying $(\mathcal{B})$ and $\mathrm{dens}(\X)<\aleph_{\omega}$. Let $\Y\subset \X$ be a closed subspace. 
Then $\mathrm{cr}(\X)\subset \mathrm{cr}(\X\quotient \Y)\cup \mathrm{dens}(\Y)^{+}$. 
In particular, such $\X$ has CSP if $\Y$ is separable and $\X\quotient \Y$ has CSP.
\begin{proof}
Observe that $\mathrm{dens}(\Y),\mathrm{dens}(\X\quotient \Y)\in\aleph_{\omega}$ are regular cardinals.
Denote the quotient map $q\colon \X\rightarrow \X\quotient \Y;\ q\colon x\mapsto x+\Y$.

Consider a regular cardinal $\kappa>\mathrm{dens}(\Y)$ such that $\kappa\in \mathrm{cr}(\X)$. 
Let $\{E_{\sigma}\}_{\sigma<\kappa}$ be the corresponding strictly nested family of 
closed subspaces of $\X$ with trivial intersection.
By Lemma \ref{cap} $(i)$ we obtain that $\bigcap_{\sigma<\kappa}\overline{\span}(E_{\sigma}\cup \Y)=\Y$.
Thus $\bigcap_{\sigma<\kappa}(\overline{\span}(E_{\sigma}\cup \Y)+\Y)=\Y$. Hence Fact \ref{quotientfact} yields that
\[\bigcap_{\sigma<\kappa}\overline{q(E_{\sigma})}=\bigcap_{\sigma<\kappa}q(\overline{\span}(E_{\sigma}\cup \Y))=0\in\X\quotient\Y.\]

Note that $E_{\sigma}\not\subset \Y$ for any $\sigma<\kappa$, since $\{E_{\sigma}\}_{\sigma<\kappa}$ is strictly nested,
$\kappa$ is regular and $\mathrm{dens}(\Y)<\kappa$.
Hence $\{0\}\subsetneq \overline{q(E_{\sigma})}$ in $\X\quotient \Y$ for each $\sigma<\kappa$. Thus by passing to a cofinal subsequence 
$\{\sigma_{\alpha}\}_{\alpha<\kappa}\subset\kappa$ such that $\{\overline{q(E_{\sigma_{\alpha}})}\}_{\alpha<\kappa}\subset \X\quotient\Y$
is strictly nested we obtain that $\kappa\in\mathrm{cr}(\X\quotient\Y)$. 

For the latter claim recall that according to Theorem \ref{subspaceseq} it holds that $\mathrm{cr}(\X)=\{\omega\}$ 
if and only if $\X$ has CSP.
\end{proof}
\end{theorem}
 
\subsection{The Kunen-Shelah properties}

In \cite{Gr} the Kunen-Shelah properties 
\[\mathrm{KS}_{7}\implies \mathrm{KS}_{6}\implies ...\implies \mathrm{KS}_{0}\]
are discussed. To mention the most important Kunen-Shelah properties in this context, a space $\X$ is said to have 
$\mathrm{KS}_{4}$, $\mathrm{KS}_{2}$ or $\mathrm{KS}_{1}$, if $\X$ admits \emph{no} uncountable polyhedron, no uncountable 
biorthogonal system or no uncountable M-basic sequence, respectively. 
It follows easily from \cite[p.114-119]{Gr} that a Banach space $\X$ with $\mathrm{KS}_{4}$ has CSP.
For example, Kunen and Shelah have provided samples of non-separable $\mathrm{KS}_{4}$ spaces by assuming CH or 
$\diamondsuit(\aleph_{1})$, see \cite[p.1086-1099]{Ne} and \cite{Sh}.
Since CSP passes to subspaces, we obtain that CSP$\implies \mathrm{KS}_{1}$ by Proposition \ref{remarkSEMI}.

It follows from recent results of Todor\v{c}evi\'{c} (see \cite{To}) that it is consistent with ZFC that the properties
$\mathrm{KS}_{7}-\mathrm{KS}_{2}$ are in fact equivalent to separability. In particular, $\mathrm{KS}_{2}\implies$CSP is consistent.
 
\begin{theorem}\label{uncountable bios}
Let $\X$ be a Banach space such that each quotient $\X\quotient \Y$ satisfies $(\mathcal{B})$ and $\mathrm{dens}(\X)<\aleph_{\omega}$. 
If $\X$ has $\mathrm{KS}_{2}$, then it has CSP.
\end{theorem}

We will first prove the following result, which will be applied.

\begin{proposition}\label{Gamma}
Let $\X$ satisfy $(\mathcal{B})$ and $\mathrm{KS}_{2}$. Suppose that $\{x_{\alpha}\}_{\alpha<\omega_{1}}\subset\X\setminus \{0\}$
and $\{\Gamma_{\sigma}\}_{\sigma<\omega_{1}}$ is a family of cofinal subsets of $\omega_{1}$
such that $\Gamma_{\alpha}\supset \Gamma_{\beta}$ for $\alpha<\beta<\omega_{1}$. Then
\[\bigcap_{\sigma<\omega_{1}}\overline{\span}(\{x_{\alpha}|\ \alpha\in\Gamma_{\sigma}\})\neq\{0\}.\]
\end{proposition}
\begin{proof}
Assume to the contrary that above
\begin{equation}\label{eq: counter}
\bigcap_{\sigma<\omega_{1}}\overline{\span}(\{x_{\alpha}|\ \alpha\in\Gamma_{\sigma}\})=\{0\}. 
\end{equation}
Then according to Lemma \ref{cap} $(ii)$
for each $\beta<\omega_{1}$ there is $\sigma<\omega_{1}$ such that 
\begin{equation}\label{eq: <cap>}
\overline{\span}(\{x_{\alpha}|\alpha<\beta\})\cap \overline{\span}(\{x_{\alpha}|\alpha\in \Gamma_{\sigma}\})=\{0\}.
\end{equation}
Let $\sigma(\beta)$ be the least ordinal satisfying \eqref{eq: <cap>} for $\sigma=\sigma(\beta)$.

Next we will define recursively an uncountable subfamily of $\{x_{\alpha}\}_{\alpha<\omega_{1}}$ by using the above notations.
Let $\alpha_{0}=0$ and for each $\theta<\omega_{1}$ let 
\[\alpha_{\theta}=\min\{\gamma\in\Gamma_{\sigma(\sup_{\epsilon<\theta}\alpha_{\epsilon})}:\ \gamma>\sup_{\epsilon<\theta}\alpha_{\epsilon}\}.\]
Note that $\{\alpha_{\theta}\}_{\theta<\omega_{1}}$ is an increasing sequence by its construction.
Observe that the corresponding family $\{y_{\theta}\}_{\theta<\omega_{1}}=\{x_{\alpha_{\theta}}\}_{\theta<\omega_{1}}$ satisfies
\[\overline{\span}(\{y_{\theta}|\theta<\gamma\})\cap\overline{\span}(\{y_{\theta}|\gamma\leq\theta\})=\{0\}\]
for each $\gamma<\omega_{1}$. 

The assumption \eqref{eq: counter} yields that 
$\bigcap_{\gamma<\omega_{1}}\overline{\span}(\{y_{\theta}|\gamma<\theta<\omega_{1}\})=\{0\}$.
Thus, an application of Lemma \ref{cap} $(i)$ for $\Y=\overline{\span}(\{y_{\theta}: \theta <\gamma\})$ 
and $Z_{\gamma}=\overline{\span}(\{y_{\theta}|\gamma<\theta<\omega_{1}\})$ for countable $\gamma$ yields the following
fact: For all $\gamma<\delta<\omega_{1}$ one can find the least ordinal 
$\eta(\gamma,\delta)\in (\delta,\omega_{1})$ such that 
\[y_{\delta}\notin\overline{\span}(\{y_{\theta}|\ \theta\in [0,\gamma]\cup [\eta(\gamma,\delta),\omega_{1})\}).\]

Put $\zeta_{0}=0$, $\zeta_{1}=1$ and for each $\alpha\in (1,\omega_{1})$ we define recursively
\[\zeta_{\alpha}=\sup_{\beta<\alpha}\eta(\sup_{\gamma<\beta}\zeta_{\gamma},\zeta_{\beta})+1.\] 
This defines an increasing sequence $\{\zeta_{\alpha}\}_{\alpha<\omega_{1}}\subset\omega_{1}$ such that
$y_{\zeta_{\gamma}}\notin \overline{\span}(\{y_{\zeta_{\theta}}|\gamma\neq \theta<\omega_{1}\})$ for $\gamma<\omega_{1}$. 

In particular $\{y_{\zeta_{\alpha}}\}_{\alpha<\omega_{1}}$ is a minimal system.
One can select by an application of the Hahn-Banach theorem suitable functionals $g_{\alpha}\in \X^{\ast}$ to 
obtain a biorthogonal system $\{(y_{\zeta_{\alpha}},g_{\alpha})\}_{\alpha<\omega_{1}}$.
This contradicts $\mathrm{KS}_{2}$, so that we have obtained the claim.
\end{proof}

\begin{proof}[Proof of Theorem \ref{uncountable bios}]
Suppose that $\X$ fails CSP and that $\mathcal{F}\subset\S_{X^{\ast}}$ is a separating set without any countable separating
subset. Then by using the proof of Theorem \ref{subspaceseq} we obtain the following:
There exists a family $\{E_{\sigma}\}_{\sigma<\lambda}$ of closed subspaces such that
$E_{\beta}\subsetneq E_{\alpha}$ whenever $\alpha<\beta<\lambda$. Here $\lambda$ is uncountable but unlike in the proof of Theorem 
\ref{subspaceseq}, here we do not need any control over the intersection $\bigcap_{\sigma}E_{\sigma}$ or the cofinality of $\lambda$. 
Moreover, similarly as in the proof of Theorem \ref{subspaceseq}, we may choose $\{E_{\sigma}\}_{\sigma}$ in the following manner. Namely, 
that there exists a family $\{f_{\phi(\gamma)}\}_{\gamma<\omega_{1}}\subset \mathcal{F}$ such that 
$E_{\sigma}=\bigcap_{\gamma<\sigma}\ker f_{\phi(\gamma)}$ 
for each $\sigma\leq \omega_{1}$. We may assign for each $\sigma<\omega_{1}$ such $x_{\sigma}\in E_{\sigma}\setminus E_{\sigma+1}$ that
$E_{\sigma+1}+[x_{\sigma}]=E_{\sigma}$ and $f_{\sigma}(x_{\sigma})=1$.

By the selection of $\{E_{\sigma}\}_{\sigma}$ and canonical identifications we get that
\[E_{\omega_{1}}=\bigcap_{\gamma<\omega_{1}}\ker f_{\phi(\sigma)}\
\mathrm{and}\ f_{\phi(\gamma)}\in (\X\quotient E_{\omega_{1}})^{\ast}=E_{\omega_{1}}^{\bot},\quad \mathrm{for}\ \gamma<\omega_{1}.\] 
Put 
\[F=\bigcap_{\beta<\omega_{1}}\overline{\span}(\{x_{\sigma}|\beta<\sigma<\omega_{1}\})\]
and observe that $F\subset E_{\omega_{1}}$, since 
$\overline{\span}(\{x_{\sigma}|\beta<\sigma<\omega_{1}\})\subset \bigcap_{\sigma<\beta}\ker(f_{\phi(\sigma)})$ for $\beta<\omega_{1}$.
Write $\hat{x}_{\sigma}=x_{\sigma}+F\in \X\quotient F$ for $\sigma<\omega_{1}$. Note that 
\begin{equation}\label{eq: 4.1}
\bigcap_{\beta<\omega_{1}}\overline{\span}(\{\hat{x}_{\sigma}|\beta<\sigma<\omega_{1}\})=\{0\}\subset \X\quotient F.
\end{equation} 

Since $\X\quotient F$ satisfies $(\mathcal{B})$ by the assumptions, we may apply the proof of Proposition \ref{Gamma}
with $\Gamma_{\sigma}=[\sigma,\omega_{1}]$ for $\sigma<\omega_{1}$ to obtain that $\X\quotient F$ admits
a biorthogonal system of length $\omega_{1}$. This can be lifted to obtain a corresponding biortohogonal system in $\X$.
\end{proof}

Regarding the assumption $\mathrm{dens}(\X)<\aleph_{\omega}$, note that if $\mathrm{dens}(\X)>2^{\omega}$, then 
$\omega^{\ast}\mathrm{-dens}(\X^{\ast})>\omega$ and $\X$ fails $\mathrm{KS}_{2}$ (see e.g. \cite[p.97-98]{Gr}). 
 
\section{Conclusions: examples, remarks and renormings}

Let us briefly recall the list of CSP spaces mentioned here: Separable spaces, the Johnson-Lindenstrauss spaces $JL_{0}$ and $JL_{2}$,
the duals $JT^{\ast}$, $JF^{\ast}$ due to Lindenstrauss-Stegal, Shelah's space $S$ under $\diamondsuit(\aleph_{1})$ and $C(K)$ spaces, 
where $K$ is the double-arrow space, Kunen's compact under CH, or any scattered separable countably tight compact.

Observe that the non-separable spaces above do not admit a system \eqref{eq: SEMIDENS}.

\begin{proposition}
Let $\X$ be a Banach space with CSP. Then $\X$ admits an equivalent uniformly Gateaux (UG) norm if and only if $\X$ is separable.
\begin{proof}
Each separable space $\X$ admits an equivalent UG norm (see \cite[Cor.6.3]{Z}). Conversely, each space $\X$ with an equivalent UG norm
is weakly countably determined and thus WLD (see \cite[Thm.6.5,Thm.3.8]{Z}). As WLD spaces are Plichko (\cite[Thm.1]{Kalenda2000}) 
we obtain by Example \ref{Plichko} that $\X$ is in fact separable.
\end{proof}   
\end{proposition}

A similar result does not hold for LUR renormings even for dual spaces. The space $JT^{\ast}$ is a non-separable CSP space, 
which admits an equivalent LUR norm by the three space property of LUR renormings (see \cite[p.1758,1785]{Z}).
On the other hand, for Kunen's compact $K$ the space $C(K)$ does not admit a Kadets-Klee and in particular not a LUR norm 
(see \cite[p.1794]{Z}).

Often the dual spaces behave better than their underlying spaces, so it is reasonable to restate the following known question.
\begin{problem}
Does any dual space $\X^{\ast}$ with CSP, i.e. a space $\X^{\ast}$ such that $\X$ is separable and does not contain $\ell^{1}$ 
isomorphically, admit an equivalent (not necessarily dual) LUR norm?
\end{problem} 
A positive answer to the previous question would partly generalize the following result.     
If $\X$ is an Asplund space then $\X^{\ast}$ has an equivalent (not necessarily dual) LUR norm (see \cite[Thm. 7.13]{Z}).
Note that any non-separable CSP dual space is \emph{not} LUR by $\omega^{\ast}$-Kadets-Klee property of $\X^{\ast}$.

If $\X^{\ast}$ has the RNP and CSP, then $\X$ is a separable Asplund space by Theorem \ref{dualchar} and
the duality of the Radon-Nikodym and Asplund properties. This means that $\X^{\ast}$ must be separable. 
\begin{problem}
Are all spaces $\X$ with CSP and the RNP in fact separable?
\end{problem} 
 
We would like to emphasize the significance of the following problems.
\begin{problem}
If $\X$ and $\Y$ have CSP, does it follow that $\X\oplus\Y$ has CSP? If so, is CSP a three-space property, i.e.
does $\X$ have CSP whenever $X\quotient \Y$ and $\Y\subset \X$ have CSP?
\end{problem}
\begin{problem}
Does there exist a CSP space without property $(C)$?
\end{problem}

\begin{proposition}\label{pd3sp}
If $\X^{\ast}$ and $\Y^{\ast}$ are dual spaces with CSP, then $\X^{\ast}\oplus \Y^{\ast}$ has CSP.
\end{proposition}
\begin{proof}
First we note that the property of being a predual of a CSP space is a three-space property.
Indeed, the property of simultaneous separability and non-containment of $\ell^{1}$ isomorphically is a three-space property 
(see \cite[2.4.h,3.2.d]{CG}). Thus we may apply Theorem \ref{dualchar}.

Now, let $\X$ and $\Y$ be some preduals of $\X^{\ast}$ and $\Y^{\ast}$, respectively. Since being a predual of a 
CSP is a three space property, we get that $\X\oplus \Y$ is a predual of a CSP space. 
The dual $(\X\oplus \Y)^{\ast}$ is thus a CSP space, which is isomorphically $\X^{\ast}\oplus \Y^{\ast}$.
\end{proof}
\subsection*{Acknowledgements}
I am grateful to Heikki Junnila for invaluable advice.

\end{document}